\documentclass[twoside,leqno,twocolumn]{article}

\usepackage[letterpaper]{geometry}

\usepackage{siamproceedings}

\usepackage{tikzscale}
\usepackage{pgfplots}
\usepackage{tikz}
\pgfplotsset{compat=newest}
\usetikzlibrary{backgrounds}
\usetikzlibrary{intersections}
\usetikzlibrary{external}
\usetikzlibrary{math}
\usetikzlibrary{positioning}
\usepackage[utf8]{inputenc}
\usepackage[T1]{fontenc}
\usepackage{amsfonts}
\usepackage{graphicx}
\usepackage{epstopdf}
\usepackage{enumitem}
\usepackage{minted}
\usepackage{algorithmic}
\usepackage{booktabs}

\ifpdf
  \DeclareGraphicsExtensions{.eps,.pdf,.png,.jpg}
\else
  \DeclareGraphicsExtensions{.eps}
\fi

\newsiamremark{remark}{Remark}
\newsiamremark{hypothesis}{Hypothesis}
\crefname{hypothesis}{Hypothesis}{Hypotheses}
\newsiamthm{claim}{Claim}

\usepackage{amsopn}

\begin{document}

\newcommand\relatedversion{}

\newif\ifblind
\blindfalse

\ifblind
  \pdfinfo{/Author ()
           /Title (Modeling and Optimization of Control Problems on GPUs)
           /Keywords (optimal control, GPU acceleration, sparse automatic differentiation, interior-point methods, Julia, nonlinear programming, domain-specific language)}
\else
  \pdfinfo{/Author (Alexis Montoison and Jean-Baptiste Caillau)
            /Title (Modeling and Optimization of Control Problems on GPUs)
           /Keywords (optimal control, GPU acceleration, sparse automatic differentiation, interior-point methods, Julia, nonlinear programming, domain-specific language)}
\fi

\ifblind
    \title{Modeling and Optimization of Control Problems on GPUs}
\else
    \title{Modeling and Optimization of Control Problems on GPUs
    \thanks{%
    Supported by the FACCTS grant ``Detecting Sparsity Patterns in Tapenade for Optimal Quantum Control Applications'' of the France-Chicago program.
    The second author is also funded by a France 2030 support managed by the Agence Nationale de la Recherche, under the reference ANR-23-PEIA-0004 (PDE-AI project).}}
\fi

\ifblind
  \author{}
\else
  \author{%
    Alexis Montoison%
    \thanks{%
      Mathematics and Computer Science Division, Argonne National Laboratory, IL, USA.
      E-mail: \email{amontoison@anl.gov}
    }
    \and
    Jean-Baptiste Caillau%
    \thanks{%
      Universit\'e C\^ote d'Azur, CNRS, Inria, LJAD.
      E-mail: \email{jean-baptiste.caillau@univ-cotedazur.fr}
    }
  }
\fi

\date{}

\maketitle


\ifblind
\fancyfoot[R]{\scriptsize{Copyright \textcopyright\ 20XX by SIAM\\
Unauthorized reproduction of this article is prohibited}}
\fi





\begin{abstract}
    We present a fully Julia-based, GPU-accelerated workflow for solving large-scale sparse nonlinear optimal control problems.
    Continuous-time dynamics are modeled and then discretized via direct transcription with \texttt{OptimalControl.jl} into structured sparse nonlinear programs.
    These programs are compiled into GPU kernels using \texttt{ExaModels.jl}, leveraging SIMD parallelism for fast evaluation of objectives, constraints, gradients, Jacobians and Hessians.
    The resulting sparse problems are solved entirely on GPU using the interior-point solver \texttt{MadNLP.jl} and the GPU sparse linear solver cuDSS, yielding significant speed-ups over CPU-based approaches.
\end{abstract}




\section{Introduction}

Solving large-scale nonlinear optimal control problems is computationally demanding, especially with fine discretizations or real-time requirements.  
While GPUs offer massive parallelism well-suited to these problems, fully exploiting their potential remains challenging due to the complexity of modeling, differentiation, and solver integration.
We present a fully GPU-accelerated workflow, entirely built in Julia~\cite{bezanson2017julia}.
Continuous-time dynamics are discretized with \texttt{OptimalControl.jl}~\cite{OC_jl} into structured, sparse nonlinear programs.  
These are compiled with \texttt{ExaModels.jl}~\cite{shin2024accelerating} into GPU kernels that preserve sparsity and compute derivatives in a single pass, enabling efficient SIMD parallelism.
Problems are solved on NVIDIA GPUs using the interior-point solver \texttt{MadNLP.jl}~\cite{shin2021graph} and the sparse linear solver \texttt{CUDSS.jl}~\cite{Montoison_CUDSS_jl_Julia_interface}, enabling end-to-end acceleration from modeling to solving.
We demonstrate the performance of this approach on benchmark problems solved on NVIDIA A100 and H100 GPUs.


\vspace{0.5cm}
\begin{tikzpicture}[
  node distance=1.3cm and 2cm,
  every node/.style={font=\scriptsize},
  box/.style={
    draw, rounded corners, thick,
    text width=2.6cm, align=center,
    minimum height=0.9cm
  },
  arrow/.style={->, thick}
]

\node[box] (oc) {OptimalControl.jl\\(continuous-time model)};
\node[box, right=of oc] (exa) {ExaModels.jl\\(GPU-friendly codegen)};
\node[box, below=of exa] (mad) {MadNLP.jl\\(GPU interior-point solver)};
\node[box, below=of oc] (cudss) {cuDSS.jl\\(GPU sparse linear solver)};

\draw[arrow] (oc) -- (exa) node[midway, above, align=center] {Discretization};
\draw[arrow] (oc) -- (exa) node[midway, below, align=center] {Sparse NLP\\ formulation};
\draw[arrow] (exa) -- (mad) node[midway, right, align=center] {SIMD-aware\\kernels};
\draw[arrow] (mad) -- (cudss) node[midway, below, align=center] {Sparse\\KKT systems};
\end{tikzpicture}

\section{Background and limitations}

Optimal control problems (OCPs) aim to find control inputs for dynamical systems modeled by ODEs that optimize a given performance criterion.
Direct transcription methods discretize these infinite-dimensional problems into large-scale nonlinear programs (NLPs).
These NLPs exhibit a sparse structure arising from time discretization: each node introduces state and control variables linked by nonlinear equality constraints enforcing the system dynamics.
Second-order methods, such as interior-point solvers, exploit this structure. 
Most existing optimal control toolchains target CPU execution.
For example, CasADi~\cite{Andersson2019} constructs symbolic expressions evaluated just-in-time or exported as C code, typically solved by CPU solvers like IPOPT~\cite{wachter2006implementation} or KNITRO~\cite{byrd2006k}, which rely on CPU linear solvers such as PARDISO~\cite{schenk2004solving}, MUMPS~\cite{amestoy2000mumps}, or HSL~\cite{fowkes2024libhsl}.
Other frameworks, such as ACADO~\cite{houska2011acado} and \texttt{InfiniteOpt.jl}~\cite{pulsipher2022unifying}, which cleverly leverage the modeling power of JuMP~\cite{dunning2017jump}, also follow the same CPU-centric paradigm.
This CPU focus limits scalability and real-time performance for large or time-critical problems that could benefit from GPU parallelism.
While some libraries provide GPU-accelerated components, none deliver a fully integrated, GPU-native workflow for nonlinear optimal control.
Our work fills this gap with a GPU-first toolchain that unifies modeling, differentiation, and solver execution, addressing the challenges of solving large-scale sparse NLPs.

\section{SIMD parallelism in direct optimal control} \label{s3}
When discretized by \emph{direct transcription}, optimal control problems (OCPs) possess an inherent structure that naturally supports SIMD parallelism. 
Consider indeed an optimal control with state $x(t) \in \mathbf{R}^n$, and control $u(t) \in \mathbf{R}^m$. Assume that the dynamics is modeled by the ODE
$$ \dot{x}(t) = f(x(t), u(t)), $$
where $f : \mathbf{R}^n \times \mathbf{R}^m \to \mathbf{R}^n$ is a smooth function. Using a one-step numerical scheme to discretise this ODE on a time grid $t_0, t_1, \dots, t_N$ of size $N + 1$ results in a set of equality constraints. For instance, with a forward Euler scheme, denoting $h_i := t_{i+1} - t_i$, one has ($X_i \simeq x(t_i)$, $U_i \simeq u(t_i)$)
$$ X_{i+1} - X_i - h_i f(X_i, U_i) = 0,\quad i = 0, \dots, N-1. $$
Similarly, a general Bolza cost that mixes endpoint and integral terms as in
$$ g(x(0), x(t_f)) + \int_0^{t_f} f^0(x(t), u(t))\,\mathrm{d}t \to \min $$
can be approximated by
$$ g(X_0, X_N) + \sum_{i=0}^{N-1} h_i f^0(X_i, U_i). $$
Discretising boundary or path constraints such as
$$ b\big(x(0),x(t_f)\big) \leq 0,\quad c\big(x(t), u(t)\big) \leq 0 $$
is obviously done according to
$$ b(X_0, X_N) \leq 0, \quad c(X_i, U_i) \leq 0,\quad i = 0, \dots, N-1. $$
The resulting NLP in the vector of unknowns $(X_0,\dots,X_N,U_0,\dots,U_{N-1})$
so involves only a few functions (\emph{kernels}), namely $f, f^0$, $g$, $b$ and $c$, that need to be evaluated on many state or control points, $X_i$, $U_i$.
This massive SIMD parallelism allows for very efficient GPU-enabled computation.
GPU acceleration thus facilitates real-time and large-scale optimal control computations critical, \emph{e.g.}, to robotics and autonomous systems as in \cite{pacaud2024gpu}.






\section{A Julia-based GPU optimization stack}
Julia offers a powerful and flexible environment for GPU programming, providing multiple levels of abstraction to suit different use cases.
The package \texttt{CUDA.jl}~\cite{besard2018juliagpu,besard2019prototyping} provides direct access to NVIDIA GPUs, supporting easy array-based programming as well as explicit CUDA kernel writing and launching.
For vendor-agnostic and portable GPU development, \texttt{KernelAbstractions.jl}~\cite{Churavy_KernelAbstractions_jl} allows writing GPU kernels in Julia that target multiple GPU vendors (NVIDIA, AMD, Intel, Apple).
This ecosystem leverages the LLVM compiler~\cite{lattner2004llvm} and vendor APIs to generate efficient native GPU code directly from pure high-level Julia code.
It allows users to exploit GPUs without requiring any knowledge of GPU programming.
For instance, \texttt{ExaModels.jl} builds on \texttt{KernelAbstractions.jl} to automatically generate specialized GPU kernels for parallel evaluation of ODE residuals, Jacobians, and Hessians needed in optimal control problems.
We build on this ecosystem to create a complete GPU-accelerated toolchain spanning modeling, differentiation, and solving.
This results into a fully Julia-native workflow for modeling and solving ODE-constrained optimal control problems on GPU.
Key components of our stack include:

\begin{itemize}
    \item[--] \texttt{OptimalControl.jl}: a domain-specific language for symbolic specification of OCPs, supporting both direct and indirect formulations.
    \item[--] \texttt{ExaModels.jl}: takes the discretized OCPs and produces sparse, SIMD-aware representations that preserve parallelism across grid points, compiling model expressions and their derivatives into optimized CPU/GPU code.
    \item[--] \texttt{MadNLP.jl}: a nonlinear programming solver implementing a filter line-search interior-point method, with GPU-accelerated linear algebra support.
    \item[--] \texttt{CUDSS.jl}: a Julia wrapper around NVIDIA’s \texttt{cuDSS} sparse solver, enabling GPU-based sparse matrix factorizations essential for interior-point methods.
\end{itemize}
\noindent Together, these components form a high-level, performant stack that compiles intuitive Julia OCP models into efficient GPU code, achieving substantial speed-ups while maintaining usability.

Our approach offers several advantages:
\begin{itemize}
    \item[--] \textbf{Abstraction}: problems can be defined intuitively using Julia DSLs (Domain Specific Languages) without any requirement for GPU programming expertise.
    \item[--] \textbf{Performance}: just-in-time compilation, SIMD parallelism, and GPU-accelerated sparse linear algebra provide substantial runtime improvements.
    \item[--] \textbf{Portability}: symbolic modeling and kernel generation are backend-agnostic; the current limitation lies in sparse linear solvers, which are still CUDA-specific, but the framework is designed to integrate alternative backends as they become available.
\end{itemize}

\section{From optimal control models to SIMD abstraction}
To illustrate the transcription from the infinite dimensional setting towards a discretized optimization suited for SIMD parallelism, consider the following elementary optimal control problem with a state function, $x(t)$, valued in $\mathbf{R}^2$, and a scalar control, $u(t)$: minimize the (squared) $L^2$-norm of the control over the fixed time interval $[0,1]$,
$$ \tfrac{1}{2}\int_0^1 u^2(t)\,\mathrm{d}t \to \min, $$
under the dynamical constraint
$$ \dot{x}_1(t) = x_2(t),\quad \dot{x}_2(t) = u(t), $$
and boundary conditions
$$ x(0) = (-1, 0),\quad x(1) = (0,0). $$
The strength of the DSL of the package \texttt{OptimalControl.jl} is to offer a syntax as close as possible to this mathematical formulation.\footnote{Note that one can actually use unicode characters to denote derivatives, integral, \emph{etc.}, making this closeness even more striking. Check \texttt{OptimalControl.jl} documentation online.} The translation of this optimal control problem so reads:

\begin{minted}{julia}
ocp = @def begin
    t in [0, 1], time
    x in R^2, state
    u in R, control
    x(0) == [-1, 0]
    x(1) == [0, 0]
    derivative(x1)(t) == x2(t)
    derivative(x2)(t) == u(t)
    integral( 0.5u(t)^2 ) => min
end
\end{minted}

The intial and final times are fixed in this case but they could be additional unknowns (see, Appendix \ref{sa1}, where the Goddard benchmark problem is modeled with a free final time. Users can also declare additional finite-dimensional parameters (or \emph{variables}) to be optimized. Furthermore, extra constraints on the state, control, or other quantities can be imposed as needed.
At this stage the crux is to seamlessly parse the abstract problem description and compile it on the fly into a discretized nonlinear optimization problem.
We achieve this by exploiting two features.
First, the DSL syntax is fully compatible with standard Julia, allowing us to use the language’s built-in lexical and syntactic parsers.
Second, pattern matching via \texttt{MLStyle.jl} \cite{MLStyle_jl} extends Julia’s syntax with additional keywords such as \verb+state+ for declaring state variables, and implements the semantic pass that generates the corresponding discretized code.
This discretized code is an \texttt{ExaModels.jl} model, which allows to declare 
optimization variables (finite dimensional vector or arrays), constraints and cost.
Regarding constraints, \texttt{ExaModels.jl} uses \emph{generators} in the form of \verb+for+ loop like statements to model the SIMD abstraction, ensuring that the function
at the heart of the statement is mapped towards a \emph{kernel} (this is where \texttt{KernelAbstractions.jl} comes into play) and efficiently evaluated by the solver. All in all, the process merely is a compilation from \texttt{OptimalControl.jl} DSL, well suited for mathematical control abstractions, into \texttt{ExaModels.jl} DSL, tailored to describe optimization problems with strong SIMD potentialities. (As explained in Section~\ref{s3}, this is indeed the case for discretizations of optimal control problems.) 
This transcription process is mostly parametrized by the numerical scheme used to discretize the ODE.
A very important outcome of having a DSL for \texttt{ExaModels.jl} models is the ability for the package to automatically differentiate the mathematical expressions involved.
Automatic differentiation (AD) is essential for modern second-order nonlinear solvers, such as IPOPT and \texttt{MadNLP.jl}, which rely on first- and second-order derivatives.

Let us take a brief look at the generated code for this simple example. The code is wrapped in a function whose parameters capture the key aspects of the transcription process: the numerical scheme (here trapezoidal), the grid size (here uniform), the backend (CPU or GPU), the initial values for variables, states, and controls (defaulting to nonzero constants across the grid), and the base precision for vectors (defaulting to 64-bit floating point):

{\small
\begin{minted}{julia}
function def(; scheme=:trapeze, grid_size=250,
  backend=CPU(), init=(0.1, 0.1, 0.1),
  base_type = Float64)
\end{minted}
}

\noindent The state declaration is compiled into an \texttt{ExaModels.jl} variable representing a $2 \times (N + 1)$ array, where $N$ is the grid size. Lower and upper bounds, plus initial values can be specified, and constraints are vectorized across grid points. Internally, the DSL uses metaprogramming to generate unique variable names and ensure proper initialization, while any syntactic or semantic errors are caught and reported at runtime.

{\small
\begin{minted}{julia}
x = begin
  local ex
  try
    variable(var"p_ocp##266", 2, 0:grid_size;
      lvar = [var"l_x##271"[var"i##275"]
      for (var"i##275", var"j##276") =
      Base.product(1:2, 0:grid_size)],
      uvar = [var"u_x##272"[var"i##275"]
      for (var"i##275", var"j##276") =
      Base.product(1:2, 0:grid_size)],
      start = init[2])
    catch ex
      println("Line ", 2, ": ",
        "(x in R^2, state)")
      throw(ex)
  end
end
\end{minted}
}


\noindent The initial-state boundary constraint must be applied across all state dimensions. This is achieved using the \verb+for+ generator.
A runtime dimension check ensures that the specified bounds match the length of the state vector being addressed:
{\small
\begin{minted}{julia}
length([-1, 0]) == length([-1, 0]) ==
  length(1:2) || throw("wrong bound dimension")
constraint(var"p_ocp##266", (x[var"i##283", 0]
  for var"i##283" = 1:2); lcon = [-1, 0],
  ucon = [-1, 0])
\end{minted}
}

\noindent The first equation of the ODE system is discretized using the trapezoidal scheme, and the corresponding expression (here the right hand side is just $x_2(t)$) is declared thanks to the \verb+for+ generator tailored for SIMD abstraction:

{\small
\begin{minted}{julia} 
constraint(var"p_ocp##266j",
  ((x[1, var"j##291" + 1] - x[1, var"j##291"])
  - (var"dt##268" * (x[2, var"j##291"] +
  x[2, var"j##291" + 1])) / 2
  for var"j##291" = 0:grid_size - 1))
\end{minted}
}

\noindent The same goes on for the second dimension of the ODE, and for the Lagrange integral cost as well, where the same numerical scheme (trapezoidal rule again) is employed for consistency (defining two objectives actually computes their sum):
 
{\small
\begin{minted}{julia}
objective(var"p_ocp##266", ((1 * var"dt##268" *
  (0.5 * var"u##277"[1, var"j##299"] ^ 2)) / 2
  for var"j##299" = (0, grid_size)))
objective(var"p_ocp##266", (1 * var"dt##268" *
  (0.5 * var"u##277"[1, var"j##299"] ^ 2)
  for var"j##299" = 1:grid_size - 1))
\end{minted}
}

\noindent The generated code returns an \texttt{ExaModels.jl} model that, when instanciated with the proper backend (\emph{e.g.}, \verb+CUDABackend()+ from \texttt{CUDA.jl}) can be passed to \texttt{MadNLP.jl} to be solved on GPU.

\section{Orchestrating model derivatives and sparse GPU linear solvers}\label{sec:orchestrate}

Our workflow executes both differentiation and linear algebra entirely on the GPU.  
The interior-point solver \texttt{MadNLP.jl} serves as the central orchestrator: it calls derivative kernels generated by \texttt{ExaModels.jl}, assembles KKT systems directly in device buffers, and invokes \texttt{CUDSS.jl} for factorization and triangular solves.  
Apart from a one-time symbolic analysis on the host, primarily reordering to reduce fill-in in the factors of the LDL$^\top$ decomposition, all operations remain on the device.
Robustness is further ensured by scaling, inertia-corrected primal–dual regularization, and optional iterative refinement.  
When vendor libraries lack suitable routines, \texttt{MadNLP.jl} falls back on \texttt{KernelAbstractions.jl} to implement custom kernels.

Each interior-point iteration therefore stays entirely on the device: derivative evaluation, KKT assembly, sparse factorizations and triangular solves, and step updates.  
All of these components are orchestrated on GPU by \texttt{MadNLP.jl}.

\section{Benchmark problems} \label{s6}

We evaluate our GPU-accelerated workflow on two challenging nonlinear optimal control problems: 
\begin{itemize}
    \item[--] the Goddard rocket problem,
    \item[--] a quadrotor control problem.
\end{itemize}
These problems have already been discretized using, for example, Euler or trapezoidal schemes, and can be found in standard optimization problem collections such as COPS~\cite{bondarenko2000cops}, as well as in the Julia port \texttt{COPSBenchmark.jl}~\cite{COPS_jl}.
 Their infinite dimensional (non discretized) counterparts are also available in the more recent open collection of control problems \texttt{OptimalControlProblems.jl} \cite{OCP_jl}. The full models using \texttt{OptimalControl.jl} are detailed in the supplementary material (check Appendix \ref{sa1}).
We compare GPU and CPU performance,
to assess both raw speed-ups and the effect of sparsity on GPU d parallelism on solver performance.
For the considered setups (see Appendix \ref{sa2} for details on the hardware and on the Julia packages used), as expected performance is better on H100 than A100, with comparable GPU over CPU speed-ups on both problems.

\begin{figure}
\includegraphics[width=.45\textwidth]{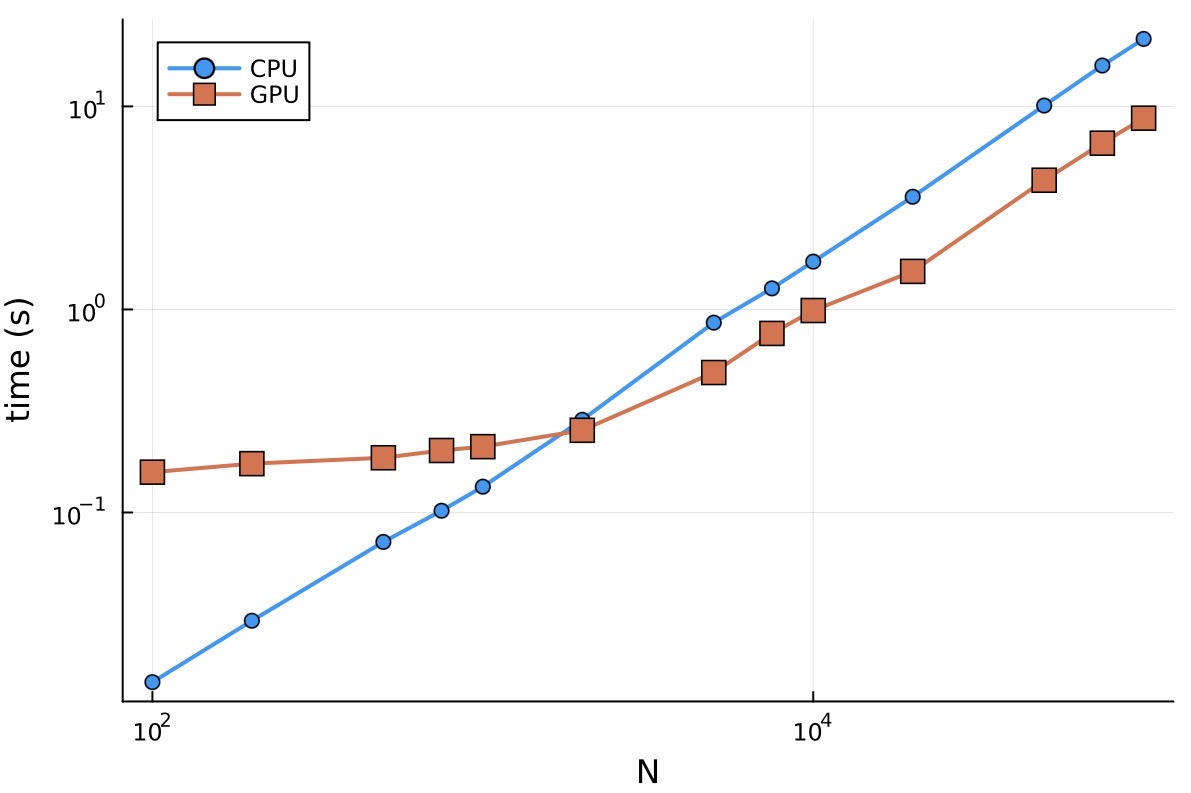}
\caption{Goddard problem, A100 solve. The grid size ranges from to $N = 1e2$ to $N = 1e5$.}
\label{fig1}
\end{figure}

On the Goddard problem, the total dimension of the discretized problem (variables plus constraints) is about $10 N$, where $N$ is the grid size. Although the problem is standard in the control literature, its solution exhibits an intricate structure for the choosen parameters with a control that is comprised of four subarcs (a concatenation of bang-singular-boundary-bang arcs). In particular, the existence of a singular arc is well known to be a source of difficulty for direct methods. Convergence towards the same objective (up to minor changes in precision depending on the number of points) is nonetheless obtained for all solves with the same (trivial) initialization, and the number of iterates remains very close on CPU and GPU for all the tested grid sizes.
On the cluster with A100, GPU is faster than CPU after $N = 2000$, see Figure~\ref{fig1}.
On the cluster with H100, GPU beats CPU after $N = 5000$, see Figure~\ref{fig2}.
On both architectures, a speed-up about $2$ is obtained.
The largest problem solved on the H100 has size about two millions with a run time about $15$~seconds.

On the Quadrotor problem, the total dimension of the discretized problem (variables plus constraints) is about $20 N$, where $N$ is the grid size.
The particular instance of the problem is unconstrained (neither control nor state
constraint, contrary to the previous Goddard example) but has a strongly nonlinear dynamics. 
As before, convergence is obtained for all architectures and grid sizes
with comparable number of iterates between CPU and GPU.
On the A100, GPU is faster than CPU after $N = 500$, see Figure~\ref{fig2}.
On the H100, GPU beats CPU after $N = 750$, see Figure~\ref{fig4}.
On both architectures, a speed-up about $5$ is obtained.
The largest problem solved on the H100 has size about $4e6$ with a run time about $13$~seconds.

 \begin{figure}
\includegraphics[width=.45\textwidth]{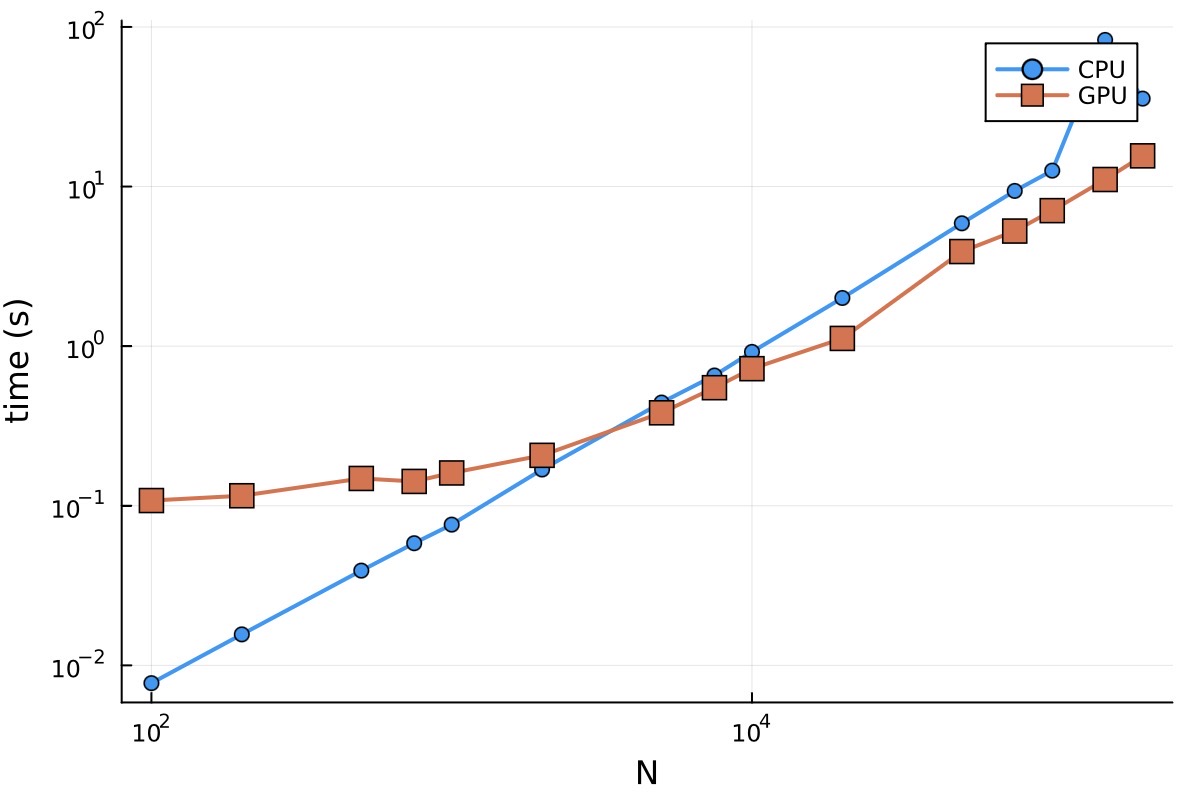}
\caption{Goddard problem, H100 solve. The grid size ranges from to $N = 1e2$ to $N = 2e5$.}
\label{fig3}
\end{figure}
 
\begin{figure}
\includegraphics[width=.45\textwidth]{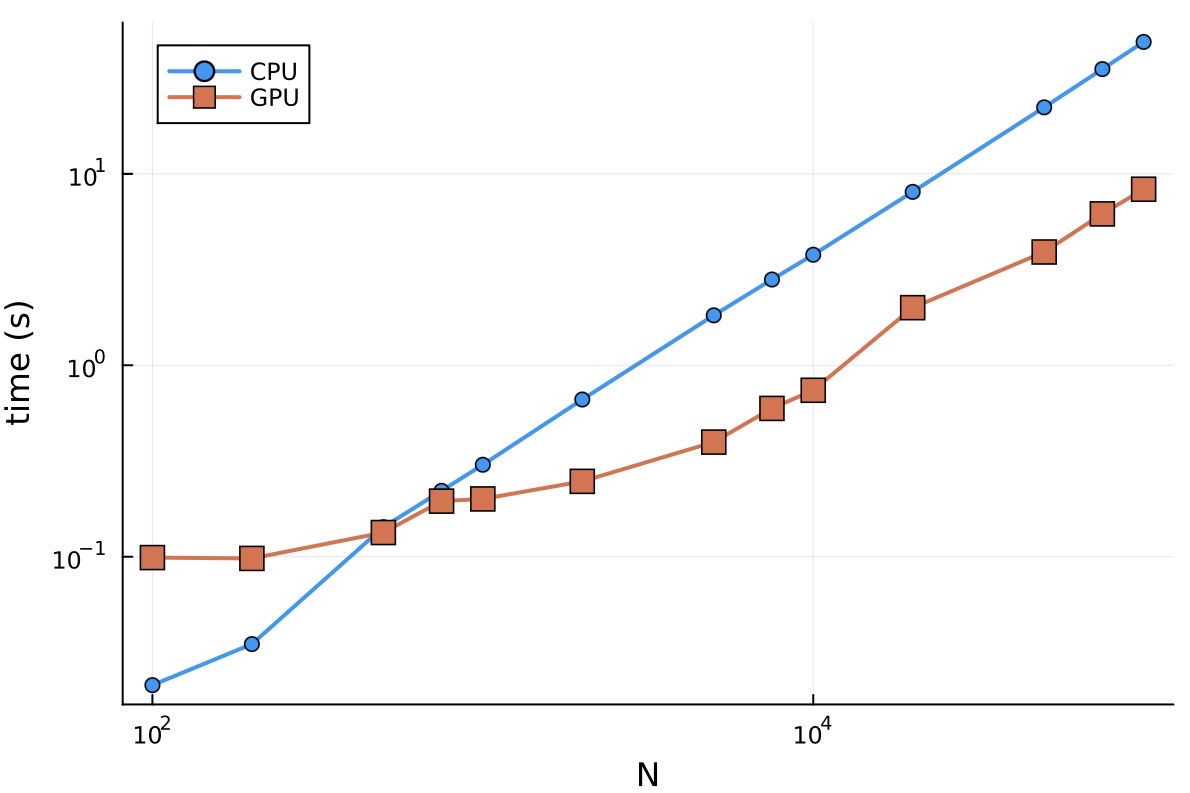}
\caption{Quadrotor problem, A100 solve. The grid size ranges from to $N = 1e2$ to $N = 1e5$.}
\label{fig2}
\end{figure}

\begin{figure}
\includegraphics[width=.45\textwidth]{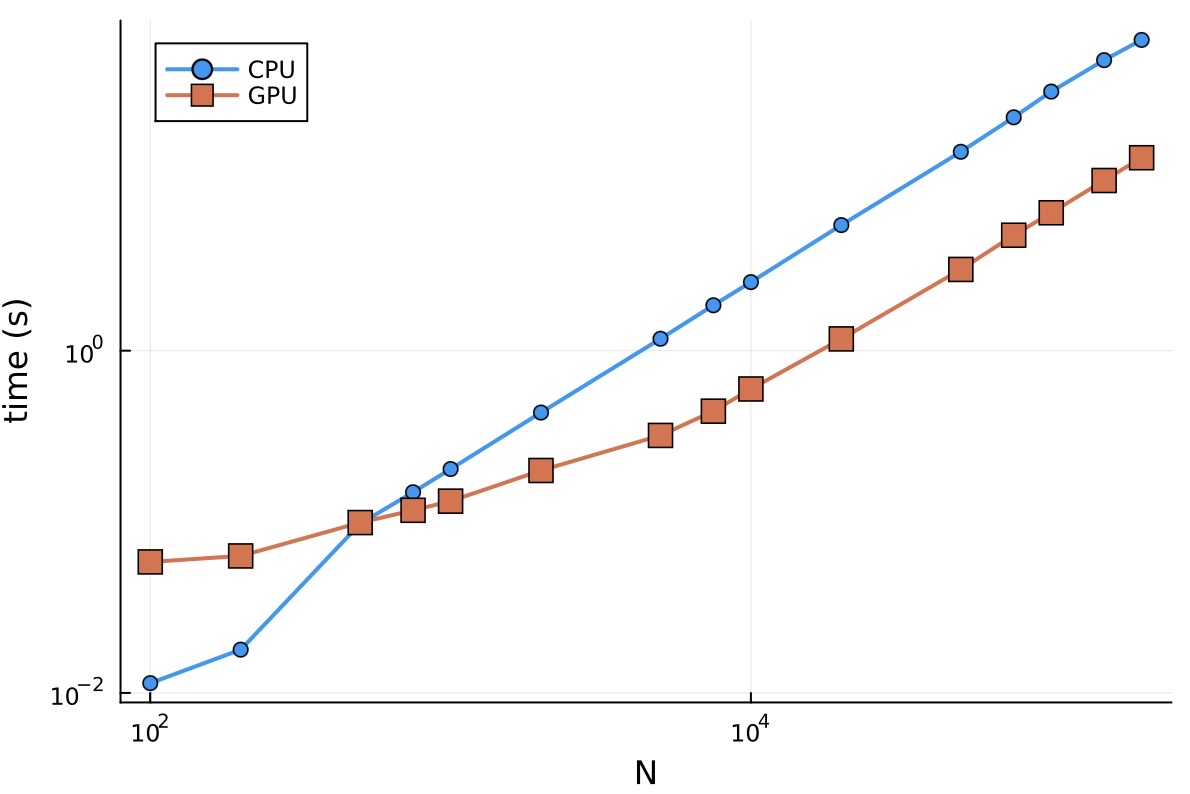}
\caption{Quadrotor problem, H100 solve. The grid size ranges from to $N = 1e2$ to $N = 2e5$.}
\label{fig4}
\end{figure}


\section{Discussion}
Julia's combination of high-level abstractions, metaprogramming, and GPU compiler infrastructure makes it uniquely suited for building performant yet expressive tools for optimal control. 
Our results show that leveraging parallelism in both model structure and solver internals unlocks substantial speed-ups, enabling new applications in aerospace engineering, quantum control, computational biology, learning, and more. 
The modular GPU stack also ensures portability across future architectures.
Future extensions include multi-GPU execution and tighter integration with differentiable programming workflows.
Overall, the synergy between Julia GPU tools and the SIMD structure of direct optimal control yields a powerful solution for solving challenging OCPs at scale.

\bibliographystyle{siamplain}
\bibliography{abbrv,main}

\begin{thebibliography}{10}

\bibitem{amestoy2000mumps}
{\sc P.~R. Amestoy, I.~S. Duff, J.-Y. L’Excellent, and J.~Koster}, {\em
  {MUMPS: a general purpose distributed memory sparse solver}}, in
  International Workshop on Applied Parallel Computing, Springer, 2000,
  pp.~121--130.

\bibitem{Andersson2019}
{\sc J.~A.~E. Andersson, J.~Gillis, G.~Horn, J.~B. Rawlings, and M.~Diehl},
  {\em {CasADi} -- {A} software framework for nonlinear optimization and
  optimal control}, Mathematical Programming Computation, 11 (2019), pp.~1--36,
  \url{https://doi.org/10.1007/s12532-018-0139-4}.

\bibitem{besard2019prototyping}
{\sc T.~Besard, V.~Churavy, A.~Edelman, and B.~De~Sutter}, {\em Rapid software
  prototyping for heterogeneous and distributed platforms}, Advances in
  Engineering Software, 132 (2019), pp.~29--46.

\bibitem{besard2018juliagpu}
{\sc T.~Besard, C.~Foket, and B.~De~Sutter}, {\em Effective extensible
  programming: Unleashing {Julia} on {GPUs}}, IEEE Transactions on Parallel and
  Distributed Systems,  (2018),
  \url{https://doi.org/10.1109/TPDS.2018.2872064},
  \url{https://arxiv.org/abs/1712.03112}.

\bibitem{bezanson2017julia}
{\sc J.~Bezanson, A.~Edelman, S.~Karpinski, and V.~B. Shah}, {\em Julia: A
  fresh approach to numerical computing}, SIAM review, 59 (2017), pp.~65--98.

\bibitem{bondarenko2000cops}
{\sc A.~S. Bondarenko, D.~M. Bortz, and J.~Mor{\'e}}, {\em {COPS: Large-scale
  nonlinearly constrained optimization problems}}, tech. report, Argonne
  National Lab., IL (US), 2000.

\bibitem{byrd2006k}
{\sc R.~H. Byrd, J.~Nocedal, and R.~A. Waltz}, {\em Knitro: An integrated
  package for nonlinear optimization}, Large-scale nonlinear optimization,
  (2006), pp.~35--59.

\bibitem{OCP_jl}
{\sc J.-B. Caillau, O.~Cots, J.~Gergaud, and P.~Martinon}, {\em
  {OptimalControlProblems.jl: a collection of optimal control problems with
  ODE's in Julia}}, \url{https://doi.org/10.5281/zenodo.17013180}.

\bibitem{OC_jl}
{\sc J.-B. Caillau, O.~Cots, J.~Gergaud, P.~Martinon, and S.~Sed}, {\em
  {OptimalControl.jl: a Julia package to model and solve optimal control
  problems with ODE's}}, \url{https://doi.org/10.5281/zenodo.13336563}.

\bibitem{Churavy_KernelAbstractions_jl}
{\sc V.~Churavy}, {\em {KernelAbstractions.jl}},
  \url{github.com/JuliaGPU/KernelAbstractions.jl}.

\bibitem{dunning2017jump}
{\sc I.~Dunning, J.~Huchette, and M.~Lubin}, {\em {JuMP: A modeling language
  for mathematical optimization}}, SIAM review, 59 (2017), pp.~295--320.

\bibitem{fowkes2024libhsl}
{\sc J.~Fowkes, A.~Lister, A.~Montoison, and D.~Orban}, {\em {LibHSL: the
  ultimate collection for large-scale scien-tific computation}}, Les Cahiers du
  GERAD ISSN, 711 (2024), p.~2440.

\bibitem{houska2011acado}
{\sc B.~Houska, H.~J. Ferreau, and M.~Diehl}, {\em {ACADO toolkit—An
  open-source framework for automatic control and dynamic optimization}},
  Optimal control applications and methods, 32 (2011), pp.~298--312.

\bibitem{lattner2004llvm}
{\sc C.~Lattner and V.~Adve}, {\em {LLVM: A compilation framework for lifelong
  program analysis \& transformation}}, in International symposium on code
  generation and optimization, 2004. CGO 2004., IEEE, 2004, pp.~75--86.

\bibitem{Montoison_CUDSS_jl_Julia_interface}
{\sc A.~Montoison}, {\em {CUDSS.jl: Julia interface for NVIDIA cuDSS}},
  \url{github.com/exanauts/CUDSS.jl}.

\bibitem{pacaud2024gpu}
{\sc F.~Pacaud and S.~Shin}, {\em {GPU}-accelerated nonlinear model predictive
  control with {ExaModels} and {MadNLP}}, arXiv e-prints,  (2024),
  pp.~arXiv--2403.

\bibitem{COPS_jl}
{\sc F.~o. Pacaud and T.~Migot}, {\em {COPSBenchmark.jl}},
  \url{github.com/MadNLP/COPSBenchmark.jl}.

\bibitem{pulsipher2022unifying}
{\sc J.~L. Pulsipher, W.~Zhang, T.~J. Hongisto, and V.~M. Zavala}, {\em A
  unifying modeling abstraction for infinite-dimensional optimization},
  Computers and Chemical Engineering, 156 (2022),
  \url{https://doi.org/doi.org/10.1016/j.compchemeng.2021.107567}.

\bibitem{schenk2004solving}
{\sc O.~Schenk and K.~G{\"a}rtner}, {\em Solving unsymmetric sparse systems of
  linear equations with pardiso}, Future Generation Computer Systems, 20
  (2004), pp.~475--487.

\bibitem{shin2024accelerating}
{\sc S.~Shin, M.~Anitescu, and F.~Pacaud}, {\em Accelerating optimal power flow
  with {GPU}s: {SIMD} abstraction of nonlinear programs and condensed-space
  interior-point methods}, Electric Power Systems Research, 236 (2024),
  p.~110651.

\bibitem{shin2021graph}
{\sc S.~Shin, C.~Coffrin, K.~Sundar, and V.~M. Zavala}, {\em Graph-based
  modeling and decomposition of energy infrastructures}, IFAC-PapersOnLine, 54
  (2021), pp.~693--698.

\bibitem{wachter2006implementation}
{\sc A.~W{\"a}chter and L.~T. Biegler}, {\em On the implementation of an
  interior-point filter line-search algorithm for large-scale nonlinear
  programming}, Mathematical programming, 106 (2006), pp.~25--57.

\bibitem{MLStyle_jl}
{\sc T.~Zhao}, {\em {MLStyle.jl}}, \url{thautwarm.github.io/MLStyle.jl}.

\end{thebibliography}

\appendix 

\section{Supplementary material}

\subsection{Descriptions of the control problems used for the benchmark}
\label{sa1}
The complete code to reproduce the runs of Section~\ref{s6} can be retrieved at the following address:
\href{https://anonymous.4open.science/r/OC-GPU-PP26}{\texttt{https://anonymous.4open.science/r/OC-GPU-PP26}}

{\small
\begin{minted}{julia}
# Goddard problem

r0 = 1.0     
v0 = 0.0
m0 = 1.0 
vmax = 0.1 
mf = 0.6   
Cd = 310.0
Tmax = 3.5

beta = 500.0
b = 2.0

o = @def begin

    tf in R, variable
    t in [0, tf], time
    x = (r, v, m) in R^3, state
    u in R, control

    x(0) == [r0, v0, m0]
    m(tf) == mf
    0 <= u(t) <= 1
    r(t) >= r0
    0 <= v(t) <= vmax

    derivative(r)(t) == v(t)
    derivative(v)(t) == -Cd * v(t)^2 *
      exp(-beta * (r(t) - 1)) / m(t) - 1 /
      r(t)^2 + u(t) * Tmax / m(t)
    derivative(m)(t) == -b * Tmax * u(t)

    r(tf) => max

end
\end{minted}
}

{\small
\begin{minted}{julia}
# Quadrotor problem

T = 1
g = 9.8
r = 0.1

o = @def begin
    
    t in [0, T], time
    x in R^9, state
    u in R^4, control

    x(0) == zeros(9)

    derivative(x1)(t) == x2(t)
    derivative(x2)(t) == u1(t) * cos(x7(t)) *
      sin(x8(t)) * cos(x9(t)) + u1(t) *
      sin(x7(t)) * sin(x9(t))
    derivative(x3)(t) == x4(t)
    derivative(x4)(t) == u1(t) * cos(x7(t)) *
      sin(x8(t)) * sin(x9(t)) - u1(t) *
      sin(x7(t)) * cos(x9(t))
    derivative(x5)(t) == x6(t)
    derivative(x6)(t) == u1(t) * cos(x7(t)) *
      cos(x8(t)) - g 
    derivative(x7)(t) == u2(t) * cos(x7(t)) /
      cos(x8(t)) + u3(t) *
      sin(x7(t)) / cos(x8(t))
    derivative(x8)(t) ==-u2(t) * sin(x7(t)) +
      u3(t) * cos(x7(t))
    derivative(x9)(t) == u2(t) * cos(x7(t)) *
      tan(x8(t)) + u3(t) * sin(x7(t)) *
      tan(x8(t)) + u4(t)

    dt1 = sin(2pi * t / T)
    dt3 = 2sin(4pi * t / T)
    dt5 = 2t / T

    0.5integral( (x1(t) - dt1)^2 +
      (x3(t) - dt3)^2 + (x5(t) - dt5)^2 +
      x7(t)^2 + x8(t)^2 + x9(t)^2 + r *
      (u1(t)^2 + u2(t)^2 +
      u3(t)^2 + u4(t)^2) ) => min

end
\end{minted}
}

\subsection{GPU detailed configurations and results} \label{sa2}
All runs performed with \texttt{OptimalControl.jl v1.1.1},
\texttt{MadNLPMumps.jl v0.5.1} and
\texttt{MadNLPGPU.jl v0.7.7}.\\

\noindent \textbf{Configuration for the A100 runs}

{\small \begin{verbatim}
julia> CUDA.versioninfo()
CUDA runtime 12.9, artifact installation
CUDA driver 12.9
NVIDIA driver 575.57.8

CUDA libraries: 
- CUBLAS: 12.9.1
- CURAND: 10.3.10
- CUFFT: 11.4.1
- CUSOLVER: 11.7.5
- CUSPARSE: 12.5.10
- CUPTI: 2025.2.1 (API 28.0.0)
- NVML: 12.0.0+575.57.8

Julia packages: 
- CUDA: 5.8.2
- CUDA_Driver_jll: 0.13.1+0
- CUDA_Runtime_jll: 0.17.1+0

Toolchain:
- Julia: 1.11.6
- LLVM: 16.0.6

1 device:
  0: NVIDIA A100-PCIE-40GB
     (sm_80, 39.490 GiB / 40.000 GiB available)
\end{verbatim}}

\noindent \textbf{Configuration for the H100 runs}

{\small \begin{verbatim}
julia> CUDA.versioninfo()
CUDA toolchain: 
- runtime 12.9, artifact installation
- driver 580.65.6 for 13.0
- compiler 12.9

CUDA libraries: 
- CUBLAS: 12.9.1
- CURAND: 10.3.10
- CUFFT: 11.4.1
- CUSOLVER: 11.7.5
- CUSPARSE: 12.5.10
- CUPTI: 2025.2.1 (API 12.9.1)
- NVML: 13.0.0+580.65.6

Julia packages: 
- CUDA: 5.8.2
- CUDA_Driver_jll: 13.0.0+0
- CUDA_Compiler_jll: 0.2.0+2
- CUDA_Runtime_jll: 0.19.0+0

Toolchain:
- Julia: 1.11.6
- LLVM: 16.0.6

Preferences:
- CUDA_Runtime_jll.version: 12.9

1 device:
  0: NVIDIA H100 80GB HBM3
     (sm_90, 79.175 GiB / 79.647 GiB available)
\end{verbatim}}

\begin{table}[htbp]
\centering
\caption{Goddard problem, A100 run}
\begin{tabular}{@{}rrr@{}}
\toprule
Grid size ($N$) & CPU time (s) & GPU time (s) \\
\midrule
100     & 1.46e-2 & 0.158 \\
200     & 2.93e-2 & 0.174 \\
500     & 7.16e-2 & 0.186 \\
750     & 1.02e-1 & 0.202 \\
1,000   & 1.34e-1 & 0.211 \\
2,000   & 2.86e-1 & 0.254 \\
5,000   & 8.61e-1 & 0.489 \\
7,500   & 1.27    & 0.762 \\
10,000  & 1.72    & 0.988 \\
20,000  & 3.59    & 1.54  \\
50,000  & 1.01e1  & 4.33  \\
75,000  & 1.59e1  & 6.59  \\
100,000 & 2.15e1  & 8.75  \\
\bottomrule
\end{tabular}
\end{table}

\begin{table}[htbp]
\centering
\caption{Quadrotor problem, A100 run}
\begin{tabular}{@{}rrr@{}}
\toprule
Grid size ($N$) & CPU time (s) & GPU time (s) \\
\midrule
100     & 2.135e-2  & 0.099    \\
200     & 3.498e-2  & 0.098    \\
500     & 1.427e-1  & 0.134    \\
750     & 2.210e-1  & 0.195    \\
1,000   & 3.023e-1  & 0.200    \\
2,000   & 6.624e-1  & 0.248    \\
5,000   & 1.823     & 0.397    \\
7,500   & 2.804     & 0.595    \\
10,000  & 3.777     & 0.740    \\
20,000  & 8.044     & 1.999    \\
50,000  & 2.223e1   & 3.904    \\
75,000  & 3.522e1   & 6.173    \\
100,000 & 4.882e1   & 8.304    \\
\bottomrule
\end{tabular}
\end{table}

\begin{table}[htbp]
\centering
\caption{Goddard problem, H100 run}
\begin{tabular}{@{}rrr@{}}
\toprule
Grid size ($N$) & CPU time (s) & GPU time (s) \\
\midrule
100     & 7.738e-3  & 0.108    \\
200     & 1.564e-2  & 0.116    \\
500     & 3.924e-2  & 0.148    \\
750     & 5.826e-2  & 0.142    \\
1,000   & 7.619e-2  & 0.161    \\
2,000   & 1.687e-1  & 0.207   \\
5,000   & 4.433e-1  & 0.382   \\
7,500   & 6.548e-1  & 0.549   \\
10,000  & 9.205e-1  & 0.723   \\
20,000  & 2.004     & 1.118   \\
50,000  & 5.884     & 3.915   \\
75,000  & 9.390     & 5.262   \\
100,000 & 1.258e1   & 7.065   \\
150,000 & 8.293e1   & 11.069  \\
200,000 & 3.556e1   & 15.575  \\
\bottomrule
\end{tabular}
\end{table}

\begin{table}[htbp]
\centering
\caption{Quadrotor problem, H100 run}
\begin{tabular}{@{}rrr@{}}
\toprule
Grid size ($N$) & CPU time (s) & GPU time (s) \\
\midrule
100     & 1.141e-2  & 0.058    \\
200     & 1.790e-2  & 0.063    \\
500     & 9.748e-2  & 0.099    \\
750     & 1.488e-1  & 0.117    \\
1,000   & 2.034e-1  & 0.132    \\
2,000   & 4.350e-1  & 0.199    \\
5,000   & 1.174     & 0.320    \\
7,500   & 1.842     & 0.443    \\
10,000  & 2.516     & 0.597    \\
20,000  & 5.420     & 1.171    \\
50,000  & 1.453e1   & 2.983    \\
75,000  & 2.311e1   & 4.728    \\
100,000 & 3.267e1   & 6.396    \\
150,000 & 4.986e1   & 9.873    \\
200,000 & 6.546e1   & 13.413   \\
\bottomrule
\end{tabular}
\end{table}

\end{document}